\documentclass[12pt, dvips, twoside]{amsart}
\usepackage{amssymb,latexsym,bbm,graphicx,epsfig,epic,eepic,oldgerm}
\usepackage{psfrag}

\theoremstyle{plain}
\newtheorem{theorem}{Theorem}[section]

\newtheorem{corollary}[theorem]{Corollary}

\newcommand{\proofend}{\hspace*{\fill} $\Box$\\}

\def\s{\smallskip}
\def\m{\medskip}
\def\eps{\varepsilon}

\def\Ham{\operatorname{Ham}}

\def\HZ{\operatorname{HZ}}

\def\gg{\gamma}
\def\gd{\delta}
\def\eps{\epsilon}

\def\gf{\varphi}

\def\go{\omega}
\def\gs{\sigma}

\def\ch{{\mathcal H}}

\def\cl{{\mathcal L}}

\def\cp{{\mathcal P}}

\def\CC{\mathbbm{C}}

\def\RR{\mathbbm{R}}

\def\pp{\partial}

\def\ra{\rightarrow}

\def\ni{\noindent}
\def\b{\bigskip}
\def\m{\medskip}

\begin{document}

\begin{titlepage}
\title{Applications of Hofer's geometry to Hamiltonian dynamics}
  
$   $ \\
$   $ \\

\address{(U.\ Frauenfelder) Department of Mathematics, Hokkaido University,
Sapporo 060-0810, Japan}
\email{urs@math.sci.hokudai.ac.jp}
\author{Urs Frauenfelder and Felix Schlenk}
\address{(F.\ Schlenk) ETH Z\"urich, CH-8092 Z\"urich, Switzerland}
\email{schlenk@math.ethz.ch}

\date{\today}

\end{titlepage}

\begin{abstract}
We prove the following three results in Hamiltonian dynamics.

\m
\begin{itemize}
\item[$\bullet$]
The Weinstein conjecture holds true for every displaceable 
hypersurface of contact type.

\s
\item[$\bullet$]
Every magnetic flow on a closed Riemannian manifold has contractible
closed orbits for a dense set of small energies.

\s
\item[$\bullet$]
Every closed Lagrangian submanifold 
whose fundamental group injects and which admits a Riemannian metric
without closed geodesics has the intersection property.
\end{itemize}

\m
\ni
The proofs all rely on the following creation mechanism for closed
orbits:
If the ray $\left\{ \gf_H^t\right\}$, $t \ge 0$, of Hamiltonian
diffeomorphisms generated by a sufficiently nice compactly supported 
time-independent Hamiltonian stops to be a minimal geodesic in its 
homotopy class,
then a non-constant contractible closed orbit must appear.
\end{abstract}
    
\maketitle

\section{Results}

\ni
We consider an arbitrary smooth symplectic manifold $(M,\go)$.
A {\it hypersurface} $S$ in $M$
is a smooth compact connected orientable codimension
$1$ submanifold of $M \setminus \pp M$ without boundary.
A closed characteristic on $S$ is an embedded circle in $S$ all
of whose tangent lines belong to the distinguished line bundle 
\[
\cl_S \,=\, \left\{ (x, \xi) \in TS \mid \go(\xi, \eta) =0 
\text{ for all } \eta \in T_x S \right\} .
\]
Examples show that $\cl_S$ might not carry any closed characteristic, 
see \cite{Gi, GG}.
We therefore follow \cite{HZ1} and consider parametrized neighbourhoods of $S$.
Since $S$ is orientable, there exists
an open neighbourhood $I$ of $0$ and a smooth diffeomorphism 
\[
\psi \colon S \times I \,\ra\, U \subset M
\]
such that $\psi (x,0) =x$ for $x \in S$.
We call $\psi$ a {\it thickening of $S$}, and we
abbreviate $S_\eps = \psi \left( S \times \left\{ \eps \right\} \right)$.
We denote by $\cp^\circ \left( S_\eps \right)$ the set of closed 
characteristics on $S_\eps$ which are contractible in $M$.
Let $\ch_c (M)$ be the set of smooth functions $[0,1] \times M \ra
\RR$ whose support is compact and contained in
$[0,1] \times \left( M \setminus \pp M \right)$.
The Hamiltonian vector field of $H \in \ch_c (M)$ defined by
\[
\go \left( X_{H_t}, \cdot \right) \,=\, -d H_t \left( \cdot \right)
\]
generates a flow $\gf_H^t$. 
The set of time-$1$-maps $\gf_H$ form the group
\[
\Ham_c (M, \go) \,:=\, \left\{ \gf_H \mid H \in \ch_c (M) \right\}
\]
of compactly supported Hamiltonian diffeomorphisms of $(M, \go)$.
We call a subset $A$ of $(M, \go)$ {\it displaceable}\,
if there exists $\gf \in \Ham_c (M,\go)$ such that 
$\gf (A) \cap A = \emptyset$.
Our basic result is

\begin{theorem}  \label{t:1}
Assume that $S$ is a displaceable hypersurface in a symplectic manifold 
$(M,\go)$.
Then for any displaceable thickening
$\psi \colon S \times I \ra U \subset M$ the set
$\left\{ \eps \in I \mid \cp^\circ \left( S_\eps \right) 
\neq \emptyset \right\}$
is dense in $I$.
\end{theorem}

Theorem~\ref{t:1} should be compared with the Hofer--Zehnder Theorem
\cite{HZ2,HZ,MaS} stating that for any thickening
$\psi \colon S \times I \ra U \subset M$ for which the Hofer--Zehnder
capacity $c_{\HZ}(U)$ is finite
the hypersurface $S_\eps$ carries a closed characteristic for almost all
$\eps \in I$.
We shall draw three applications from Theorem~\ref{t:1}.

\b
\ni
{\bf 1. The Weinstein conjecture}

\s
\ni
A hypersurface $S$ in $(M, \go)$ is called {\it of contact type}\, 
if there exists a Liouville vector field $X$ 
(i.e., $\cl_X \go = d \iota_X \go = \go$)
which is defined in a neighbourhood of $S$
and is everywhere transverse to $S$. 
Weinstein conjectured in \cite{W} that every hypersurface $S$ 
of contact type with $H^1(S;\RR) =0$ carries a closed characteristic. 

\begin{corollary}  \label{c:wein}
Assume that $S \subset (M,\go)$ is a displaceable hypersurface of
contact type. Then $\cp^\circ (S) \neq \emptyset$.
In particular,
the Weinstein conjecture holds true for every displaceable hypersurface
of contact type.
\end{corollary}


Proofs of the Weinstein conjecture for {\it all}\, hypersurfaces of
contact type of
{\it special classes}\, of symplectic manifolds have been found in
\cite{V, HZ1, HV, FHV, HV2, J, Ma, V1, Lu, V2, C, LT, Lu2}.
Corollary~\ref{c:wein} generalizes the results in \cite{V, HZ1, FHV,
V2, LT}, where the ambient symplectic manifold is of the form 
$\left( M \times \RR^2, \go \oplus \go_0 \right)$.
Under the additional assumption that 
$(M, \go)$ satisfies $[\go] |_{\pi_2(M)} =0$ and is convex, 
Corollary~\ref{c:wein} has been proved in \cite{FS}.

\b
\ni
{\bf 2. Closed trajectories of a charge in a magnetic field}

\s
\ni
Consider a Riemannian manifold $(N,g)$ of dimension at least $2$,
and let $\go_0 = \sum_i dp_i \wedge dq_i$ be the standard symplectic form on 
the cotangent bundle $T^*N$.
The motion of a unit charge on $(N,g)$ subject
to a magnetic field can be
described as the Hamiltonian flow of the Hamiltonian system
\begin{equation}  \label{e:Ham}
H \colon \left( T^*N, \go_\gs \right) \ra \RR, 
\quad\, H (q,p) \mapsto \frac12 \left| p \right|^2,
\end{equation}
where $\gs$ is the closed $2$-form on $N$ representing the magnetic field
and where $\go_\gs$ is the twisted symplectic form
$\go_\gs = \go_0 + \pi^* \gs$.
A trajectory of a charge on $(N,g)$ in the magnetic field $\gs$ has
constant speed, and closed trajectories $\gg$ on $N$ of speed $c >0$
correspond to closed orbits of \eqref{e:Ham}
on the energy level $E_c = \left\{ H = c^2 /2 \right\}$.
An old problem in Hamiltonian mechanics asks for closed orbits 
on a given energy level $E_c$, see \cite{Gi0}.
We denote by $\cp^\circ \left( E_c \right)$ the set of closed
orbits on $E_c$ which are contractible in $T^*N$ and hence
project to contractible closed trajectories on $N$.

\begin{corollary}  \label{c:mag}
Consider a closed Riemannian manifold $(N,g)$ endowed with a  
closed $2$-form $\gs$ which does not vanish identically. 
There exists $d>0$ such that $\cp^\circ \left( E_c \right) \neq
\emptyset$ for a dense set of $c \in \;]0, d[$.
\end{corollary}

The number $d>0$ has a geometric meaning: If the Euler characteristic
$\chi(N)$ vanishes, $d$ is the supremum of the real numbers $c$ for
which the sublevel set
\[
H^c \,=\, \left\{ (q,p) \in T^*N \mid H(q,p) = \tfrac 12 |p|^2 \le c
\right\} 
\]
is displaceable in $\left( T^*N, \go_\gs \right)$, and if $\chi(N)$ does
not vanish, $d$ is defined via stabilizing \eqref{e:Ham} by 
$\left( T^* S^1, dy \wedge dx \right) \ra \RR$, $(x,y) \mapsto \frac 12 |y|^2$.
Corollary~\ref{c:mag}
generalizes various other results on closed orbits 
of magnetic flows on small energy levels:
The existence of a sequence $c \ra 0$ with 
$\cp^\circ \left( E_c \right) \neq \emptyset$
has been proved by Polterovich \cite{P0} and Macarini \cite{Mac} 
under the assumption that $[\gs] |_{\pi_2(M)} =0$
and by Ginzburg--Kerman \cite{GK2} under the assumption that $\gs$ is
symplectic. 
If $\gs$ is symplectic and $[\gs] |_{\pi_2(M)} =0$,
Corollary~\ref{c:mag} has been proved in \cite{CGK}.
For exact forms \cite{FS} and rational symplectic forms 
\cite{GG,Mac,Mac2} it is known that
$\cp^\circ (E_c) \neq \emptyset$ for almost all sufficiently small $c>0$,
and if $\gs$ is symplectic and compatible with $g$, then 
$\cp^\circ (E_c) \neq \emptyset$ for all sufficiently small $c>0$, see
\cite{K}.
We refer to \cite{Gi0} and \cite[Section 12.4]{FS}
for more details on the state of the art of the existence problem 
for closed trajectories of a charge in a magnetic field.

\b
\ni
{\bf 3. Lagrangian intersections}

\s
\ni
A middle-dimensional submanifold $L$ of $(M, \go)$ is called {\it Lagrangian}\,
if $\go$ vanishes on $L$. 
According to a celebrated theorem of Gromov, \cite[2.3.$B_3'$]{Gr}, 
a closed Lagrangian submanifold $L \subset M \setminus \pp M$
with $[\go] |_{\pi_2(M,L)}=0$ is not displaceable.
The following result generalizes Theorem~1.4.A in \cite{LP} 
and Theorem~13.1 in \cite{FS}.

\begin{corollary}  \label{c:lag}
Assume that $L \subset M \setminus \pp M$ is a closed Lagrangian 
submanifold of $(M, \go)$ such that
\begin{itemize}
\item[(i)]
the injection $L \subset M$ induces an injection $\pi_1(L) \subset
\pi_1(M)$;  
\s
\item[(ii)]
$L$ admits a Riemannian metric none of whose closed geodesics is contractible.
\end{itemize}
Then $L$ is not displaceable.
\end{corollary}

\b
\ni
{\bf Acknowledgements.}
We cordially thank Leonid Polterovich for pointing out to us that the use
of functions as in the figure below in his approach to closed
orbits of a charge in a magnetic field will provide a result like 
Corollary~\ref{c:mag}.
This paper was written during the second authors stay at Tel Aviv University
in April 2003. He wishes to thank Hari and Harald and TAU for their warm 
hospitality.

\section{Proofs}

\ni
{\it Proof of Theorem~\ref{t:1}:}
We closely follow Polterovich's beautiful argument in \cite[Section 9.A]{P0}.
We fix $\eps_0 \in I$ and choose $\gd >0$ so small that 
$[\eps_0-2\gd, \eps_0+2\gd] \subset I$.
Let $f \colon \RR \ra [0,1]$ be a smooth function such that
\[
f(\eps) =0 \,\text{ if } \eps \notin \;]\eps_0-2\gd, \eps_0+2\gd[ ,
\quad
f(\eps) = 1 \,\text{ if } \eps \in [\eps_0-\gd, \eps_0+\gd] ,
\]
and
\[
f'(\eps) > 0 \,\text{ if } \eps \in \; ]\eps_0-2\gd, \eps_0-\gd[ ,
\quad
f'(\eps) <0 , \,\text{ if } \eps \in \; ]\eps_0+\gd, \eps_0+2\gd[ ,
\]
see the figure below.

\begin{figure}[h] 
 \begin{center}
   \psfrag{ee}{$\eps$}
   \psfrag{e}{$\eps_0$}
   \psfrag{-2ed}{$\eps_0-2\gd$}
   \psfrag{-ed}{$\eps_0-\gd$}
   \psfrag{ed}{$\eps_0+\gd$}
   \psfrag{2ed}{$\eps_0+2\gd$}
   \psfrag{f}{$f(\eps)$}
   \psfrag{1}{$1$}
   \leavevmode\epsfbox{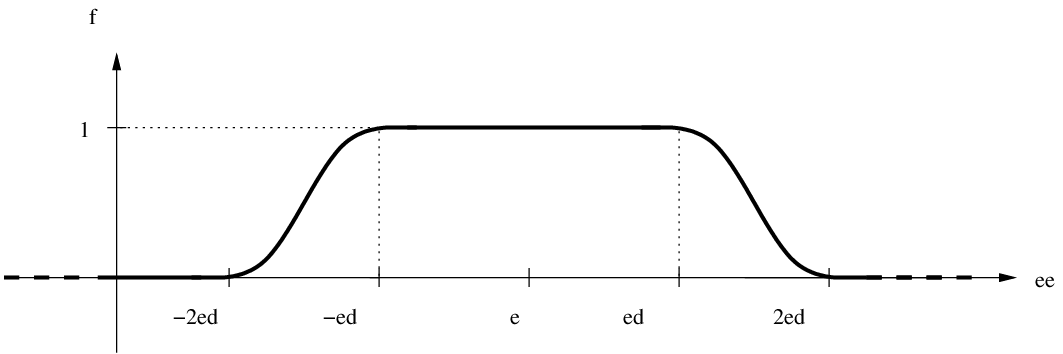}
 \end{center}
\end{figure}
%
%
%

\noindent
We define the time-independent Hamiltonian $F \in \ch_c \left( M
\right)$ by
\begin{equation*}  
F (x) \,=\, 
 \left\{
  \begin{array}{ll}
   f(\eps) & \text{if } x \in S_\eps , \\ [0.2em]
   0    & \text{otherwise}.
  \end{array}
 \right.
\end{equation*}
The norm $\| H \|$ of $H \in \ch_c(M)$ is defined as
\[
\| H \| \,=\, \int_0^1 
\left( \sup_{x \in M} H(t,x) - \inf_{x \in M} H(t,x) \right) dt .
\]
We say that $H,K \in \ch_c(M)$ are equivalent, $H \sim K$, 
if $\gf_H = \gf_K$ and the paths $\left\{ \gf_H^t \right\}$, 
$\left\{ \gf_K^t \right\}$, $t \in [0,1]$, are homotopic with fixed end points.
Following \cite{P0} we set
\[
\mu \left( F,s \right) \,=\, 
\inf \left\{ \left\| H \right\| \mid H \sim sF \right\},
\quad\, s \ge 0 ,
\]
and define the asymptotic non-minimality $\mu(F)$ of $F$ as
\[
\mu (F) \,=\, 
\lim_{s \ra \infty} \frac{\mu \left( F,s\right)}{\left\| F \right\|} .
\]
Then $\mu (F) \in [0,1]$.
Since the support of $F$ is displaceable, it follows from
\cite[Theorem 3.3.A]{BP} or \cite[Theorem 8.3.A]{P} that
 \[
\mu(F) \,\le\, \frac 12 ,
\]
see also \cite{LM1}.
On the other hand, our choice of $F$ guarantees that the Hessian of $F$ at any
of its critical points vanishes.
If the flow $\gf_F^s$, $s \ge 0$, had no non-constant contractible 
closed orbit, then $sF$ were ``slow'' for every $s \ge 0$ in the terminology 
of \cite{MSl}, and so, according to Theorem~1.4 in \cite{MSl},
\[
\mu (F) \,=\, 1 ,
\]
see also \cite{E,Mc, Oh0}.
This contradiction shows that the flow $\gf_F^s$ has a non-constant 
contractible closed orbit,
and so there exists $\eps \in \;]\eps_0 -2\gd, \eps_0 +2\gd[$
such that $\cp^0 \left( S_\eps \right) \neq \emptyset$.
Since $\eps_0 \in I$ and $\gd >0$ sufficiently small were arbitrary,
Theorem~\ref{t:1} follows.
\proofend

\ni
{\it Proof of Corollary~\ref{c:wein}:} 
Choose $\eps >0$ so small that the local flow $\psi_t$ of the Liouville vector
field $X$ defined near $S$ exists for all $t \in\; ]-\eps, \eps[$.
Then $\psi_t$ induces a bijection 
$\cp^\circ \left( S \right) \ra \cp^\circ \left( \psi_t (S) \right)$,
$x \mapsto \psi_t x$, see \cite[p.\ 122]{HZ}.
Corollary~\ref{c:wein} thus follows from Theorem~\ref{t:1}.
\proofend

\ni
{\it Proof of Corollary~\ref{c:mag}:}
We closely follow \cite[Section 12.1]{FS}. 
Let $(N,g)$ and $\left(T^*N, \go_{\gs} \right)$ be as in 
Corollary~\ref{c:mag}.
Since $\gs$ does not vanish, $\dim N \ge 2$, and so every energy level
$E_c = \left\{ H = c^2/2 \right\}$, $c>0$,
is a connected hypersurface.
Let $\chi (N)$ be the Euler characteristic of $N$.

\s
\ni
{\bf Case 1.} $\chi (N) =0$.
We define $d \in [0,\infty[\:\cup\;\! \{ \infty \}$ by
\begin{equation*}  
d \,=\, 
d(g, \gs) \,=\, \sup \left\{ c \ge 0 \mid H^c \text{ is displaceable in }
\left( T^*N, \go_\gs \right) \right\} .
\end{equation*}
Notice that 
\[
d \,=\, \sup \left\{ c \ge 0 \mid E_c \text{ is displaceable in }
\left( T^*N, \go_\gs \right) \right\} .
\]
Since $\gs \neq 0$, the zero section $N$ of $T^*N$ is not Lagrangian,
and so a remarkable theorem of Polterovich \cite{P2,LS} implies that $d>0$.
We can assume without loss of generality that $d$ is 
finite.\renewcommand{\thefootnote}{\fnsymbol{footnote}}\footnote[2]{If 
$\gs$ is exact, then $d(g,\gs)$ is not larger than Man\'e's critical value
and in particular finite, see \cite[Section 12.2]{FS}.
If $\gs$ is the area form on the flat $2$-torus $T^2$, then
$\left( T^*T^2, \go_\gs \right)$ is symplectomorphic to 
$\left( T^2, \gs \right) \times \left( \RR^2, \go_0 \right)$, 
see \cite[Section 2.4.B]{EP},
and so $d(g,\gs) =\infty$.} 
\renewcommand{\thefootnote}{\arabic{footnote}}Corollary~\ref{c:mag} 
follows from applying Theorem~\ref{t:1} to 
$S = E_{d/2}$ and a thickening
\[
\psi \colon S \times \left] -d/2,d/2 \right[  \,\ra\, \bigcup_{0<c<d} E_c
\]
such that $\psi \left( S \times \{ \eps \} \right) = E_{\eps + d/2}$.

\s
\ni
{\bf Case 2.} $\chi (N) \neq 0$.
In this case the zero section $N$ is not displaceable for topological
reasons. 
We use a stabilization trick used before by Macarini \cite{Mac}.
Let $S^1$ be the unit circle, and  denote canonical coordinates on
$T^*S^1$ by $(x,y)$.  
We consider the manifold 
$T^* \left( N \times S^1 \right) = T^* N \times T^* S^1$ 
endowed with the split symplectic form 
$\go = \go_\sigma \oplus \go_{S^1}$, 
where $\go_{S^1} = dy \wedge dx$.
Then $N \times S^1$ is not Lagrangian, and 
$\chi \left( N \times S^1 \right) =0$.
Let
\[
H_1(q,p) = \tfrac{1}{2} |p|^2,
\quad
H_2(x,y) = \tfrac{1}{2} |y|^2,
\quad
H(q,p,x,y) =  \tfrac 12 |p|^2 + \tfrac 12 |y|^2 
\]
be the metric Hamiltonians on $T^*N$, $T^*S^1$ and $T^*N \times T^*S^1$.
In order to avoid confusion, we denote their energy levels by
$E_c(H_1)$, $E_c(H_2)$ and $E_c(H)$.
Repeating the argument given in Case~1 for the Hamiltonian system
\begin{equation}  \label{e:m:H}
H \colon \left( T^*N \times T^* S^1, \go \right) \ra \RR 
\end{equation} 
and 
\[
d \,=\, d(g, \gs) \,=\, \sup \left\{ c \ge 0 \mid H^c \text{ is
displaceable in } \left( T^*N \times T^* S^1, \go \right) \right\} 
\]
we find that 
$\cp^\circ \left( E_c(H) \right) \neq \emptyset$ for a dense set of 
$c \in \;]0,d[$. 
Fix $c \in \;]0,d[$ such that $\cp^\circ \left(E_c (H) \right)
\neq \emptyset$.
Since the Hamiltonian system \eqref{e:m:H} splits,
a contractible closed orbit $x(t)$ on $E_c (H)$ is of the form
$\left( x_1(t), x_2(t) \right)$, where $x_1$ is a contractible
closed orbit on $E_{c_1} (H_1)$ and $x_2$ is a contractible closed orbit
on $E_{c_2} (H_2)$ and $c_1 + c_2 = c$. 
Since the only contractible orbits of $H_2 \colon T^*S^1 \ra \RR$
are the constant orbits on $E_0 (H_2)$, we conclude that $c_2 = 0$ 
and $c_1 = c$, and so $x_1 \in \cp^\circ \left( E_c (H_1)\right)$.
The proof of Corollary~\ref{c:mag} is complete.
\proofend

\ni
{\it Proof of Corollary~\ref{c:lag}:} 
Arguing by contradiction we assume that $L$ and hence a neighbourhood 
$U$ of $L$ in $M$ is displaceable.
We choose a Riemannian metric as in (ii) and denote by $|p|$ the length
of a covector $(q,p) \in T^*L$.
By Weinstein's Theorem we find $\eps >0$ such that a neighbourhood
of $L$ in $U$ can be symplectically identified with 
$T_{2\eps}^*L = \left\{ (q,p) \in T^*L \mid |p| < 2 \eps \right\}$.
Then
$S_\eps = \left\{ (q,p) \in T^*_{2\eps}L \mid |p| = \eps \right\}$
is a displaceable hypersurface of contact type, and hence, by
Corollary~\ref{c:wein}, carries a closed characteristic which is 
contractible in $M$. It corresponds to
a closed orbit of the geodesic flow on 
$S_\eps$ which is contractible in $M$.
According to (i) and (ii) such orbits to not exist, however.
\proofend

\enddocument